\documentclass[11pt]{article}
\renewcommand{\baselinestretch}{1.5}

\def\sizesmallfig{0.4}
\usepackage{subfigure}
\usepackage{amsmath,amssymb}
\usepackage{graphicx}
\usepackage{makeidx}

\newenvironment{pf}{\noindent {\bf Proof. }}{\hfill $\square$}

\newcommand{\bea}{\begin{eqnarray*}}
\newcommand{\eea}{\end{eqnarray*}}

\newcommand{\bydef}{\stackrel{\bigtriangleup}{=}}

\newcommand{\vect}[1]{{\boldsymbol #1 }}

\newcommand{\inprod}[2]{\langle #1 , #2 \rangle }
\newcommand{\card}[1]{\vert #1 \rvert}
\newcommand{\bc}{\begin{center}}
\newcommand{\ec}{\end{center}}

\newcommand{\by}{\mbox{\boldmath $y$}}
\newcommand{\refs}[1]{(\ref{#1})}

\newcommand{\bx}{\vect{x}}

\newcommand{\R}{\mathbb R}
\newcommand{\N}{\mathbb N}

\newcommand{\be}{\begin{equation}}
\newcommand{\ee}{\end{equation}}
\newcommand{\beaa}{\begin{eqnarray*}}
\newcommand{\eeaa}{\end{eqnarray*}}
\newcommand{\ben}{\begin{enumerate}}
\newcommand{\een}{\end{enumerate}}
\newcommand{\db}{\hspace*{\fill}{\zapf o}}

\newcommand{\bpn}{\begin{proposition}\twlsf}
\newcommand{\epn}{\db\end{proposition}}
\newcommand{\bdm}{\begin{displaymath}}
\newcommand{\edm}{\end{displaymath}}
\newcommand{\ba}{\begin{array}}
\newcommand{\ea}{\end{array}}

\newcommand{\mb}[1]{\boldsymbol{#1}}

\newcommand{\st}{\mathop{\rm s.t.}}

\newcommand{\MP}{{\cal P}}
\newcommand{\dd}{\mathrm{d}}
\newcommand{\mbal}{\mb{\alpha}}
\newcommand{\mbbt}{\mb{\beta}}
\newcommand{\mbga}{\mb{\gamma}}
\newcommand{\mbth}{\mb{\theta}}
\newcommand{\mbmu}{\mb{\mu}}

\newtheorem{definition}{Definition}
\newtheorem{lemma}{Lemma}
\newtheorem{proposition}{Proposition}

\newtheorem{remark}{Remark}
\newtheorem{theorem}{Theorem}

\newcommand{\eps}{\epsilon}
\newcommand{\Om}{\Omega}

\usepackage{color}
\definecolor{dkviolet}{rgb}{0.6,0,0.8}

\def\jean#1{\textcolor{red}}
\newcommand{\new}[1]{\text{{\color{black}#1}}}

\title{Optimal data fitting: a moment approach\footnote{The research of the first author was funded by the European Research Council (ERC) under the European's Union Horizon 2020 research and innovation program (grant agreement 666981 TAMING). The research of the second author has been partially supported by the LabEx PERSYVAL-Lab (ANR-11-LABX-0025-01) funded by the French program ``Investissement d'avenir'' and by the European Research Council (ERC) ``STATOR'' Grant Agreement nr. 306595. }} 
\author{Jean-Bernard Lasserre \thanks{LAAS, 7 Avenue Du Colonel Roche, 31077 Toulouse C\'edex 4, France, lasserre@laas.fr}
\and Victor Magron \thanks{L2S CENTRALESUPELEC; 3, Rue Joliot-Curie,  91192 Gif sur Yvette; France, victor.magron@l2s.centralesupelec.fr}}

\date{\today}

\begin{document}
\maketitle

\begin{abstract}
We propose a  moment relaxation for two problems, the separation and covering problems with semi-algebraic sets generated by a polynomial of degree $d$. We show that (a)  the optimal value of the  relaxation finitely converges to the optimal value of the original  problem, when the moment order $r$ increases and (b) after performing some small perturbation of the original problem, convergence can be achieved with $r=d$. We further provide a practical iterative algorithm  that is computationally tractable for large datasets and present encouraging computational results. 
\end{abstract}

\section{Introduction}
\label{sec:intro}
Data fitting problems have long been very useful in many different application areas. A well-known problem is the problem of finding the minimum-volume ellipsoid in $n$-dimensional space $\R^n$ containing all points that belong to a given finite set $S \subset \R^n$. This minimum-volume covering ellipsoid problem is important in the area of robust statistics, data mining, and cluster analysis (see e.g. Sun and Freund \cite{SunFreund04} and the recent book by M. Todd \cite{Todd16}). Pattern separation as described in Calafiore \cite{Calafiore02} is another related problem, in which an ellipsoid that separates a set of points $S_1$ from another set of points $S_2$ needs to be found under some appropriate optimality criteria such as minimum volume or minimum distance error.

{
These problems have been studied for a long time. The minimum-volume covering ellipsoid (MVCE) problem was discussed by John \cite{John48} in 1948. This problem has been modeled as a convex optimization problem with linear matrix inequalities (LMI) and solved by interior-point methods (IPM) in Vandenberghe et al. \cite{Vandenberghe98}, Sun and Freund \cite{SunFreund04} and Magnani et al. \cite{Magnani05}.
In particular, the "dual-reduced-Newton" algorithm presented in~\cite{SunFreund04} combines interior-point and active-set methods allowing one to efficiently compute the MVCE of moderately large datasets (in dimension $n = 30$ and dataset with $30000$ points, it takes a few seconds on a personal laptop).  The recent survey by Todd~\cite{Todd16} provides more details on algorithms depending on the size of the datasets and the dimension. In particular, for huge-scale problems ($n = 500$ and and datasets with $500 000$ points), the Wolfe-Atwood algorithm~\cite{Todd07} is the only one able to compute the MVCE.}

\begin{figure}[!htbp]
\centering
\includegraphics[scale=\sizesmallfig]{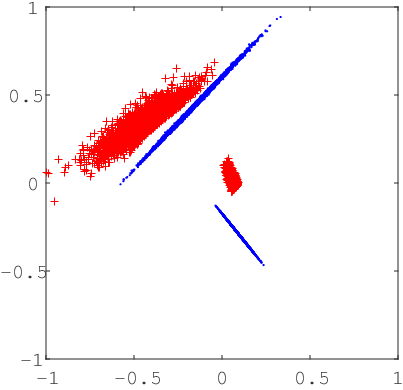}
\caption{Two datasets $S_1$ (blue) and $S_2$ (red) non separable with ellipsoids.}
\label{fig:nonsepbyellips}
\end{figure}

{The problem of pattern separation via ellipsoids was studied by Rosen \cite{Rosen65} and Barnes \cite{Barnes82}. Glineur \cite{Glineur98} has proposed some methods to solve this problem with different optimality criteria via conic programming. Although efficient algorithms are already available, they could become computationally intractable if the cardinality of datasets is large. In addition, different semi-algebraic sets other than ellipsoidal sets could be considered for these data fitting problems. Note that for complicated clusters as in Figure~\ref{fig:nonsepbyellips}, it will be impossible in general to separate two datasets with ellipsoids.}



%

{This was our motivation to approximate such datasets with the level set of a polynomial $\theta$ of a priori fixed degree $d$, possibly greater than 2. Given two datasets $S_1$ and $S_2$, the superlevel (resp.~sublevel) set of $\theta$ should contain $S_1$ (resp.~$S_2$). Interestingly, this approach can also be used for minimum-volume covering problems when $S_2 = \emptyset$. However for very large datasets, this may not be competitive with dedicated algorithms, such as the ones presented in~\cite{Todd16}.}

With $\mbga = (\gamma_1,\dots,\gamma_n) \in \N^n$, let us note $\mb{x}^{\mbga} \bydef  x_1^{\gamma_1} \, \cdots \, x_n^{\gamma_n}$ and $|\mbga| \bydef \gamma_1 + \dots + \gamma_n$.
One possible approach to find the coefficients of the above-mentioned polynomial $\theta = \sum_{\mbga} \theta_{\mbga}\mb{x}^{\mbga}$ is to solve the following linear problem (LP):

\begin{equation}
\label{eq:lp}
\left[
\ba{lrr}
\inf & \sum_{|\mbga|\leq d} |\theta_{\mbga}| & \quad\\
\st & \theta(\mb{x})\geq 0, & \forall \mb{x}\in S_1,\\
\quad & -\theta(\bx)\geq 0, & \forall \bx\in S_2,\\
\quad &  \theta \in \R[\mb{x}], & \deg \theta = d.
\ea
\right]
\end{equation}
\if{
min  sum_alpha |theta_alpha|
s.t. \forall x in S_1, \sum_\alpha \theta_alpha(x^alpha) >= 0 \,,
     \forall x in S_2, \sum_\alpha \theta_alpha(x^alpha) <= 0 \,.
}\fi    
When $\theta$ has degree $d$, this LP has $\binom{n+d}{n}$ variables and $|S_1| + |S_2|$ constraints on the vector of coefficients of the polynomial $\theta$. The variables of LP~(1) are the coefficients $(\theta_\gamma)_{|\gamma| \leq d}$ of the polynomial $\theta$. \\
Given a feasible solution $\theta$ of LP~(1), the superlevel (resp.~sublevel) set of $\theta$ contains $S_1$ (resp.~$S_2$).
Therefore, we choose the $\ell_1$-norm of the coefficient vector of $\theta$ for the objective function, in order to minimize the volume of the level-sets of $\theta$.\\
This LP may become ill-conditioned whenever several points from the clusters are close to each other. The reason is that in this case the LP has almost redundant inequality constraints.
We performed practical experiments with several LP solvers (\new{Gurobi~\cite{gurobi}}, SDPT3~\cite{sdpt3}, Mosek~\cite{mosek}, SeDuMi~\cite{sedumi}), which all include a pre-processing step to remove nearly dependent constraints (see~\cite[\S~1.3.5]{sdpt3}).

By solving LP~\eqref{eq:lp}, we were able to separate datasets of small size (less than $ 10^2$ points). For various randomly generated datasets of larger size, such as the datasets $S_1$ and $S_2$ (with $10^5$ points) depicted in Figure~\ref{fig:nonsepbyellips}, we experienced numerical issues either with the algorithm implemented in the LP solvers.\footnote{\new{For instance, the Gurobi solver cannot avoid numerical issues due to the large magnitude of matrix coefficients.} The SDPT3 solver (version 3.4.0) is not able to compute the solution of LP~\eqref{eq:lp} and often aborts with various error messages, including the following one: ``\texttt{Stop: steps are too short}''.}\\
Another drawback of this LP formulation is that it cannot tackle all data fitting applications considered in the present study, including minimum-volume covering ellipsoids.
\subsection*{Contributions and Paper Outline}
{In this paper 
we propose a common methodology for these data fitting problems and in particular its extension to 
general basic semi-algebraic sets (more general than ellipsoids) based on a moment approach.}
This methodology is based on the {\em moment-SOS approach} and has the distinguishing feature to {\em not} individualize each point of the two  clouds of data points, that is, we do not 
incorporate positivity constraints of the type $\pm\theta(\mb{x}_i)\geq0$ for each
point of the cloud. More precisely:

\begin{enumerate}
\item[(1)] In Section \ref{sec:momentformulation}, we propose a  {\em moment
 relaxation} for these data fitting problems with basic semi-algebraic sets $\Omega=\{\bx\in\R^n:\theta(\bx)\geq 0\}$, where $\theta\in\R[\bx]$ is a polynomial with a priori fixed degree $d$. One main advantage 
 is to avoid considering $|S_1| + |S_2|$ constraints (i.e. avoids individualizing each point). The information of each dataset $|S_i|$ is collected in a localizing matrix $M_r(\theta \mb{y}^i)$ associated to an empirical measure $\mu^i$ supported on $|S_i|$ (and where $y^i$ is a finite vector of moments associated with $\mu^i$). In our case, we perform a smoothing thanks to the two LMIs:
 $M_r(\theta \mb{y}^1) \succeq 0$ (associated to $S_1$) and $M_r(\theta \mb{y}^2) \preceq 0$ (associated to $S_2$), with $r \in \N$. When the first (resp.~second) condition is satisfied for all $r$, this is actually equivalent to the nonnegativity of $\theta$ (resp.~$-\theta$) on the support of $\mu^1$ (resp.~$\mu^2$) (see~\cite{Lasserre11}). We show in Proposition~\ref{prop:relaxation} that the optimal value of the  relaxation converges in finitely many steps to the optimal value of the original  problem, when the moment order $r$ increases. The key idea of our approach is that,
instead of imposing a constraint for each point
 in the dataset, we  require that the support of any
 probability measure $\mu$ that is generated on the dataset is contained in $\Omega$.
Using powerful results from the theory of moments,
we may replace all membership constraints by two Linear Matrix Inequality (LMI) constraints of size $\binom{n+r}{r}$. Hence for 3D-datasets 
the size of each of the two LMIs is $O(r^3/6)$ (and so for a quartic polynomial ($r = 4$), the size of each LMI is only $\binom{3+4}{3} = 35$). 
\item[(2)] In Section~\ref{sec:exactconvergence}, we show the following result: If 
$\max[\vert S_1\vert,\vert S_2\vert]=s$  and the degree $d$
of the polynomial $\theta$  is such that $\binom{n+d}{n}\geq s$, 
then finite convergence occurs at step $r=d$, generically.
This genericity condition can be ensured after performing some arbitrary small perturbation of the original problem. The possible drawbacks of this method is that for large size clusters, the size of the localizing matrices grows rapidly, leading to LMIs which are expensive to solve. Therefore to handle large datasets in practice we combine the above approach with 
a heuristic inspired from results on  semi-infinite optimization by Ben-Tal et al.~\cite{Ben-Tal79}.

\item[(3)]  In Section \ref{sec:algorithms}, we provide a practical iterative
algorithm based on the results of Section \ref{sec:exactconvergence} for
these data fitting problems that is computationally
tractable for datasets with a very large number of points. The corresponding method is an iterative procedure where the degree of the polynomial is fixed in advance (typically $r = 2$ or $r = 4$) and where we solve a moment relaxation with measures supported on subsets of the initial cluster. 
Even though this algorithm does not always converge, it happens to be very efficient in practice. We present encouraging computational results of cluster separation with up to $10^5$ points, either with ellipsoid or quartic level sets. 
\end{enumerate}

\section{Moment Relaxations}
\label{sec:momentformulation}
\subsection{Problem Formulation}
\label{ssec:pformulation}
With $\mb{x} = (x_1,\dots,x_n)$, consider a polynomial $\theta \in \R[\mb{x}]$ of degree at most $d$: $\theta(\mb{x})=\sum_{\mbga\in\N^n:\card{\mbga}\leq d}\theta_{\mbga}\mb{x}^{\mbga}$. 
Letting $\mbth = \{\theta_{\mbga}:\mbga\in\N^n, \card{\mbga}\leq d\}$ be the coefficient vector of $\theta$, $\mbth\in\R^t$, where $t = \binom{n+d}{d}$ and $\Theta\subset\R^t$, we obtain a family of semi-algebraic sets $\Omega_{\theta}=\{\bx\in\R^n:\theta(\bx)\geq 0\}$ for $\theta\in\Theta$. The problem of separating a finite dataset $S_1 \subset \R^n$ from another finite dataset $S_2 \subset \R^n$ with one of these semi-algebraic sets can be written as follows: 
\begin{equation}
\label{eq:sprob}
{\MP}^s \quad \left[
\ba{llrr}
\tau^s = & \inf & f(\mbth) & \quad\\
&\st & \theta(\mb{x})\geq 0, & \forall \mb{x}\in S_1,\\
&\quad  & -\theta(\bx)\geq 0, & \forall \bx\in S_2,\\
&\quad  & \mbth\in\Theta ,
\ea
\right]
\end{equation}
where $f$ is an optimality criterion and $\tau^s$ is the optimal value of $\MP^s$.

If we only consider one dataset $S$, then we can formulate the problem of covering $S$ with the best semi-algebraic set $\Omega_{\theta}$ with respect to optimality criterion $f$ as follows:
\begin{equation}
\label{eq:cprob}
{\MP}^c \quad \left[
\ba{llrr}
\tau^c = & \inf & f(\mbth) & \quad\\
& \st & \theta(\mb{x})\geq 0, & \forall \mb{x}\in S,\\
& \quad & \mbth\in\Theta .
\ea
\right]
\end{equation}

If $f$ is the volume function of ellipsoids and $\theta$ is a quadratic function that generates ellipsoidal sets, then the pattern separation via ellipsoids and minimum-volume covering ellipsoid problems are obtained respectively from these two general problems. Thus, for $r = 1$, we  consider $f = \log \det \mb{Q}^{-1}$ with $\mb{Q}$ being a positive definite matrix of size $n$ and a separating polynomial $\theta(\mb{x}) = -\mb{x}'\mb{Q}\mb{x}+\mb{b}'\mb{x} + c$. 
When using level sets of polynomials with higher degree $2 r$, we consider the same cost function $f = \log \det \mb{Q}^{-1}$ with a positive definite matrix $\mb{Q}$ of size $\binom{n+r-1}{n}$ and a separating polynomial $\theta(\mb{x}):= -v_r(\mb{x})'\mb{Q} v_r(\mb{x})+\mb{b}'v_r(\mb{x}) + c$, where $v_r(\mb{x})$ is the vector of degree-$r$ monomials, i.e.~$(x_1^r,x_1^{r-1}x_2,\dots, x_n^r)$. 
Note that we can also consider the more general separating polynomial $\theta(\mb{x}):= 1 -w_r(\mb{x})'\mb{Q} w_r(\mb{x})$, with $w_r(\mb{x})$ being the vector of all monomials with degree at most $r$ (See Section~\ref{ssec:mvce} and Section~\ref{ssec:spve} for more practical details).
%

Since the covering problem is a special case of the separation problem ($S_2=\emptyset$), we focus on the latter problem in the following sections. 

\subsection{Moment Formulation}
\label{ssec:formulation}
We now investigate the application of the moment-SOS approach (see Henrion \cite{Henrion04}, Lasserre \cite{Lasserre01},   and the references therein) to Problem (\ref{eq:sprob}). Let $\mu^i$ be a probability measure generated on $S_i$, $i=1,2$,
\begin{equation}
\label{eq:mudef}
\mu^i := \sum_{\mb{x}\in S_i} \mu^i_{\mb{x}}\delta_{\mb{x}},
\end{equation}
where $\delta_{\mb{x}}$ denotes the Dirac measure at $\mb{x}$, $\sum_{\mb{x}\in S_i}\mu^i_{\mb{x}}=1$, and $\mu^i_{\mb{x}}\geq 0$ for all $\mb{x}\in S_i$, $i=1,2$. For example, the uniform probability measure $\mu^i$ generated on $S_i$ has $\mu^i_{\mb{x}} = 1/\card {S_i}$ for all $\mb{x}\in S_i$.

All the moments $\mb{y}^i=\left\{y^i_{\mbal}\right\}$ of $\mu^i$ are calculated as follows:
\begin{equation}
\label{eq:ydef}
y^i_{\mbal} = \int \mb{x}^{\mbal} \dd \mu^i = \sum_{\mb{x}\in S_i}\mu^i_{\mb{x}}\mb{x}^{\mbal}, \quad \mbal \in \N^n.
\end{equation} 
For any nonnegative integer $r$, the $r$-moment matrix associated with $\mu^i$ (or equivalently, with $\mb{y}^i$) $M_r(\mu^i)\equiv M_r(\mb{y}^i)$ is a matrix of size $\binom{n+r}{r}$. Its rows and columns are indexed in the canonical basis $\left\{\mb{x}^{\mbal}\right\}$ of $\R [\mb{x}]$, and its elements are defined as follows:
\begin{equation}
\label{eq:mmdef}
M_r(\mb{y}^i)(\mbal,\mbbt) = y^i_{\mbal +\mbbt}, \quad \mbal, \mbbt \in \N^n, \card {\mbal}, \card {\mbbt} \leq r.
\end{equation}
Similarly, given $\theta \in \R [\mb{x}]$, the localizing matrix $M_r(\mbth \mb{y}^i)$ associated with $\mb{y}^i$ and $\theta$ is defined by
\begin{equation}
\label{eq:lmdef}
M_r(\mbth\mb{y}^i)(\mbal,\mbbt) := \sum_{\mbga \in \N^n}\theta_{\mbga}y^i_{\mbal+\mbbt+\mbga}, \quad \mbal, \mbbt \in \N^n, \card {\mbal}, \card {\mbbt} \leq r,
\end{equation}
where $\mbth = \left\{\theta_{\mbga} \right\}$ is the vector of coefficients of $\theta$ in the canonical basis $\left\{\mb{x}^{\mbal} \right\}$. 

If we define the matrix $M_r^{\mbga}(\mb{y}^i)$ with elements
$$
M_r^{\mbga}(\mb{y}^i)(\mbal,\mbbt) = y^i_{\mbal +\mbbt+\mbga}, \quad \mbal, \mbbt, \mbga \in \N^n, \card {\mbal}, \card {\mbbt} \leq r,
$$
then the localizing matrix can be expressed as $M_r(\mbth\mb{y}^i)=\sum_{\mbga\in \N^n}\theta_{\mbga}M_r^{\mbga}(\mb{y}^i)$.

Note that for every polynomial $f \in \R [\mb{x}]$ of degree at most $r$ with its vector of coefficients denoted by $\mb{f} = \left\{f_{\mbga}\right\}$, we have:
\begin{equation}
\label{eq:lmprop}
\inprod{\mb{f}}{M_r(\mbth\mb{y}^i)\mb{f}} = \int \theta f^2 \dd \mu^i .
\end{equation}
This property shows that necessarily, $M_r(\mbth\mb{y}^i)\succeq \mb{0}$, whenever $\mu^i$ has its support contained in the level set $\left\{\mb{x}\in \R^n: \theta(\mb{x})\geq 0\right\}$. Therefore, if we replace all membership constraints in Problem $\MP^s$ by two LMI constraints $M_r(\mbth\mb{y}^1)\succeq 0$ and $M_r(-\mbth\mb{y}^2)\succeq 0$, we obtain a relaxation of Problem $\MP^s$:
\begin{equation}
\label{eq:formulation}
{\MP}^s_r(\mb{y}^1,\by^2) \quad \left[
\ba {lrl}
\tau^s_r(\mb{y}^1,\by^2) = & \inf & f(\mbth)\\
& \st & M_r(\mbth\mb{y}^1)\succeq \mb{0}\\
& \quad & M_r(-\mbth\mb{y}^2)\succeq \mb{0}\\
& \quad & \mbth\in\Theta ,
\ea
\right]
\end{equation}
with optimal value denoted by $\tau^s_r(\mb{y}^1,\by^2)$.
We emphasize that ${\MP}^s_r(\mb{y}^1,\by^2)$ depends on $\mb{y}^1$ and $\mb{y}^2$, the respective moment sequences of the two measures $\mu^1$ and $\mu^2$, that are both fixed beforehand. Next, we prove that the convergence of $\tau^s_r(\mb{y}^1,\by^2)$ to $\tau^s$ occurs under  mild properties of $\mu^1$ and $\mu^2$.

\subsection{Convergence as the number of moments increases}
\label{ssec:finiteconvergence}
Compared to $\MP^s$, the data of $S_1$ and $S_2$ are aggregated into the vector
$\mb{y}^1$ and $\by^2$ used in ${\MP}^s_r(\mb{y}^1,\by^2)$. Both problems have exactly the same variables, but
\begin{enumerate}
\item[-] Problem $\MP^s$ has $\card {S_1}+\card{S_2}$ linear constraints, whereas
\item[-] Problem ${\MP}^s_r(\mb{y})$ has two LMI constraints $M_r(\mbth\mb{y}^1)\succeq 0$ and $M_r(-\mbth\mb{y}^2)\succeq 0$ with matrix size $\binom{n+r}{r}$.
\end{enumerate}
If $r$ is not too large, solving ${\MP}^s_r(\mb{y})$ is preferable to
solving $\MP^s$, especially if $\card {S_1}+\card{S_2}$ is large. It is natural to ask 
 how good this moment relaxation could be as compared to the original
problem and which value of $r$ we have to use to obtain a strong lower bound. In this section, let us assume that fixed probability measures $\mu^i$ generated on $S_i$, $i=1,2$, are selected; for example, the uniform probability measures as mentioned in the previous section.
\begin{proposition}
\label{prop:relaxation}
Let $\theta\in\R[\bx]$, and let $\MP^s$, ${\MP}^s_r(\mb{y}^1,\by^2)$, $r\in \N$ be as in \refs{eq:sprob} and \refs{eq:formulation} respectively. Then:
$$
\tau^s_r(\mb{y}^1,\by^2)\leq \tau^s_{r+1}(\mb{y}^1,\by^2),\quad \mbox{and}\quad \tau^s_r(\mb{y}^1,\by^2)\leq
\tau^s, \quad \forall r\in \N ,
$$
where $\tau^s$ and $\tau^s_r(\mb{y}^1,\by^2)$ are optimal values of $\MP^s$ and ${\MP}^s_r(\mb{y}^1,\by^2)$ respectively. 
\end{proposition}
\begin{pf}
For every $\mbga \in \N^n$, $M_r^{\mbga}(\mb{y}^i)$ is the north-west corner square sub-matrix with size $\binom{n+r}{r}$ of $M_{r+1}^{\mbga}(\mb{y}^i)$, $i=1,2$. This follows directly from the definition of the matrix $M_r^{\mbga}(\mb{y}^i)$. 

Since $M_r(\mbth\mb{y}^1)=\sum_{\mbga\in \N^n}\theta_{\mbga}M_r^{\mbga}(\mb{y}^1)$ for all $r$,  $M_r(\mbth\mb{y}^1)$ is also a north-west corner square submatrix of $M_{r+1}(\mbth\mb{y}^1)$. This implies that if $M_{r+1}(\mbth\mb{y}^1)\succeq \mb{0}$, then $M_r(\mbth\mb{y}^1)\succeq \mb{0}$. Similar arguments can be applied for $M_r(-\mbth\mb{y}^2)$ and $M_{r+1}(-\mbth\mb{y}^2)$. Thus, any feasible solution of ${\MP}^s_{r+1}(\mb{y}^1,\by^2)$ is feasible for ${\MP}^s_r(\mb{y}^1,\by^2)$. So we have:
$$
\tau^s_r(\mb{y}^1,\by^2)\leq \tau^s_{r+1}(\mb{y}^1,\by^2), \quad \forall r\in \N .
$$
Similarly, any feasible solution of $\MP^s$ is feasible for  ${\MP}^s_r(\mb{y}^1,\by^2)$. Indeed, if $\mbth$ is feasible for $\MP^s$ then we have $\theta(\mb{x})\geq 0$ for all $\mb{x}\in S_1$ and $\theta(\bx)\leq 0$ for all $\bx\in S_2$. Therefore, the probability measures $\mu^1$ and $\mu^2$ defined in \refs{eq:mudef} have their supports contained in the level set $\left\{\mb{x}\in \R^n: \theta(\mb{x})\geq 0\right\}$ and $\left\{\mb{x}\in \R^n: \theta(\mb{x})\leq 0\right\}$ respectively. In view of \refs{eq:lmprop}, we have $M_r(\mbth\mb{y}^1)\succeq \mb{0}$ and $M_r(-\mbth\mb{y}^2)\succeq \mb{0}$. This proves that $\mbth$ is feasible for ${\MP}^s_r(\mb{y}^1,\by^2)$ with any $r\in \N$. Thus,
$
\tau^s_r(\mb{y}^1,\by^2)\leq \tau^s, \quad \forall r\in \N .
$
\end{pf}

We next show that 
 if $\mu^i$ is supported on the whole set $S_i$, $i=1,2$, then the optimal values $\left\{\tau^s_r(\mb{y}^1,\by^2)\right\}$ converges to $\tau^s$, when $r$ increases and the convergence is finite. The statement is formally stated and proved as follows:
\begin{theorem}
\label{thr:finiteconvergence}
Let $\MP^s$, $\mu^i$ and ${\MP}^s_r(\mb{y}^1,\by^2)$, $r\in \N$ be as in \refs{eq:sprob}, \refs{eq:mudef} and \refs{eq:formulation}, respectively. If $\mu^i_{\mb{x}}>0$ for all $\mb{x}\in S_i$, $i=1,2$, then 
$$
\tau^s_r(\mb{y}^1,\by^2)\uparrow \tau^s\quad\mbox{ as $r$ increases,}
$$ 
and the convergence is finite.
\end{theorem}
\begin{pf}
From Proposition \ref{prop:relaxation}, we have
$
\tau^s_r(\mb{y}^1,\by^2) \leq \tau^s, \quad \forall r\in \N .
$
In addition, as $\mu^i$ in \refs{eq:mudef} is finitely supported, its moment matrix $M_r(\mb{y}^i)$ defined in~\refs{eq:mmdef} with $\mb{y}^i$ as in \refs{eq:ydef} has finite rank. That is, there exists $r^i_0\in \N$ such that 
$$
\mbox{rank}(M_r(\mb{y}^i)) = \mbox{rank}(M_{r^i_0}(\mb{y})), \quad \forall r\geq r^i_0.
$$
In other words, $M_r(\mb{y}^i)$ is a flat extension of $M_{r_0}(\mb{y}^i)$ for all $r\geq r^i_0$ (see Curto and Fialkow \cite{CurtoFialkow00} for more details).

Now, let $r_0 :=\max\{r^1_0,r^2_0\}$ and let $\mbth$ be an arbitrary $\eps$-optimal solution of ${\MP}^s_{r_0}(\mb{y}^1,\by^2)$, $\eps>0$,
i.e., $f(\mbth) \leq \tau^s_{r_0}(\mb{y}^1,\by^2) + \eps$. As $M_{r_0}(\mbth\,\mb{y}^i)\succeq \mb{0}$ and since $M_{r_0 + \frac{\deg \theta}{2}}(\mb{y}^i)\succeq \mb{0}$, we may invoke Theorem 1.6 in  \cite{CurtoFialkow00} and deduce that $\mu^1$ has its support contained in the level set $\left\{\mb{x} \in \R^n: \theta(\mb{x})\geq 0\right\}$. Similarly, $\mu^2$ has its support contained in the level set $\left\{\mb{x} \in \R^n: \theta(\mb{x})\leq 0\right\}$. This implies that $\theta(\mb{x})\geq 0$ for all $\mb{x}\in S_1$ and $\theta(\mb{x})\leq 0$ for all $\mb{x}\in S_2$ because $\mu^i$ is supported on the whole set $S_i$ ($\mu^i_{\mb{x}}>0$ for all $\mb{x}\in S_i$). Thus, $\mbth$ is feasible for $\MP^s$ and $\tau^s\leq f(\mbth) \leq \tau^s_{r_0}(\mb{y}^1,\by^2) + \eps$. We then have
$$
\tau^s_{r_0}(\mb{y}^1,\by^2)\leq \tau^s\leq \tau^s_{r_0}(\mb{y}^1,\by^2) + \eps.
$$
As $\eps>0$ was arbitrary, $\tau^s_{r_0}(\mb{y}^1,\by^2)=\tau^s$. Since from 
 Proposition \ref{prop:relaxation} $\tau^s_r(\mb{y}^1,\by^2)$ is monotone and
 bounded, we obtain that $\tau^s_r(\mb{y}^1,\by^2)\uparrow \tau^s$
and the convergence is finite.
\end{pf}

Theorem \ref{thr:finiteconvergence} provides the  rationale for solving
${\MP}^s_r(\mb{y}^1,\by^2)$ instead of $\MP^s$. However, despite the finite
convergence
we do not know how large the value of $r_0$ could be. In
the next section, we will discuss how to select appropriate values $r$
for Problem ${\MP}^s_r(\mb{y}^1,\by^2)$.
\section{Convergence of Measures}
\label{sec:exactconvergence}
In this section, we analyze how the genericity of the datasets $S_1 \subset \R^n$ and $S_2 \subset \R^n$ affects the convergence of ${\MP}^s_r(\mb{y}^1,\by^2)$.

\begin{definition}[$r$-genericity]
\label{def:rgeneric}
We say that $S$ is said ``$r$-generic'' when $S$ does not belong to the level set of any polynomial of degree at most $r$.
\end{definition}
We next investigate how different 
(and much ``smaller") atomic probability measures can be selected 
to  yield optimality.
\subsection{Convergence in Measure}
\label{ssec:exactconvergence}
As mentioned in the previous section,  $\tau^s_r(\mb{y}^1,\by^2)$ converges to
$\tau^s$ as $r$ increases, and there exists an index $r_0$ such that  $\tau^s_{r_0}(\mb{y}^1,\by^2)=\tau^s$. As for several other problems reformulated with moment relaxations, 
no explicit value for the relaxation order $r_0$ is available in general.
We investigate next the dependence of $r_0$ 
 on $\card {S_1}$ and $\card{S_2}$ and show that 
 under certain rank conditions,
 we will have $\tau^s_r(\mb{y}^1,\by^2)=\tau^s$ for any probability measures $\mu^i$ supported on the whole set $S_i$, $i=1,2$.
\begin{proposition}
\label{prop:exactconvergence}
Let $\MP^s$, $\mu^i$, and ${\MP}^s_r(\mb{y}^1,\by^2)$ be as in \refs{eq:sprob}, \refs{eq:mudef}, and \refs{eq:formulation}, respectively. Assume that $\card{S_i} = \binom{n+r}{r}$, $S_i$ is $r$-generic (in the sense of Definition~\ref{def:rgeneric}), 
and $\mu^i_{\mb{x}}>0$ for all $\mb{x}\in S_i$ for $i=1,2$. Then the following holds: 
\[
\card{S_i} = \mbox{rank} \, (M_r(\mb{y}^i)) 
\quad \text{and} \quad
\tau^s_r(\mb{y}^1,\by^2)=\tau^s \,.
\]
\end{proposition}
\begin{pf}
Note that  $M_r(\mb{y}^i)$ is an 
 $\binom{n+r}{r}\times \binom{n+r}{r}$ matrix and the probability measure
 $\mu^i$ is supported on the whole set $S_i$, with $\card {S_i} =
 \binom{n+r}{r}$, i.e.~$\mu^i$ is an $\card {S_i}$-atomic measure.
For $i = 1,2$, we first show that the rank of the matrix $M_r(\mb{y}^i)$ is  maximal, i.e.~$\mbox{rank} \, (M_r(\mb{y}^i)) = \binom{n+r}{r} = \card{S_i}$.  For each $x \in S_i$, we denote by $\zeta_x := (x^\alpha)_{|\alpha|\leq r}$ the moment sequence associated to the Dirac measure at $x$. Note that $M_r(\mb{y}^i) = \sum_{x \in S_i} \mu_x \, \zeta_x \zeta_x^T = G G^T$, where $G$ is the $\binom{n+r}{r}\times \binom{n+r}{r}$ matrix whose columns are the vectors $\sqrt{\mu_x} \zeta_x$. We show that $G$ is invertible. Indeed, if $G z = 0$ for some vector $(z_\alpha)_{|\alpha|\leq r}$, then one has $\sqrt{\mu_x} \sum_{|\alpha| \leq r} z_\alpha x^\alpha = 0$, for all $x \in S_i$. Thus, $S_i$ belongs to the level set $\{\mb{x}: g_i(\mb{x})=0\}$ of the polynomial $g_i$ (of degree at most $r$) with vector of coefficients $z$. This contradicts our assumption and therefore necessarily
$z = 0$, which implies that $G$ is non singular and which in turn
implies that $M_r(\mb{y}^i)$ is also invertible. Hence, one has $\card{S_i} = \mbox{rank} \, (M_r(\mb{y}^i))$.

Then according to Laurent \cite[Lemma~2.7]{Laurent05}, there exist $\card {S_i}$ interpolation polynomials $f_j\in \R [\mb{x}]$ of degree at most $r$, $j=1,\ldots,\card {S_i}$, such that
$$
f_j(\mb{x}(k)) = \left\{
\ba{lr}
0, & j\neq k,\\
1, & j=k,
\ea
\right. \quad \forall j,k=1,\ldots,\card {S_i} .
$$
Now, let $\mbth$ be an arbitrary $\eps$-optimal solution of ${\MP}^s_r(\mb{y}^1,\by^2)$, $\eps>0$. As $\mbth$ is feasible for ${\MP}^s_r(\mb{y}^1,\by^2)$, $M_r(\mbth \mb{y}^1)\succeq \mb{0}$, $M_r(-\mbth \mb{y}^2)\succeq \mb{0}$ and $f(\mbth) \leq \tau^s_r(\mb{y}^1,\by^2) + \eps$. For every $j = 1,\ldots,\card {S_1}$, we have:
$$
\inprod{\mb{f}_j}{M_r(\mbth \mb{y}^1)\mb{f}_j}\geq 0,
$$
where $\mb{f}_j$ is the vector of coefficients of the polynomial $f_j$. Then
Eq. \refs{eq:lmprop} implies that 
$$
\int\theta f_j^2 \dd\mu^1 \geq 0\Leftrightarrow \mu^1_{\mb{x}(j)}\theta(\mb{x}_j) \geq 0.
$$
Since $\mu^1_{\mb{x}}>0$ for all $\mb{x}\in S_1$,  we obtain $\theta(\mb{x}(j))\geq 0$ for all $j=1,\ldots,\card{S_1}$. Similarly, we also obtain $\theta(\mb{x}(j))\leq 0$ for all $j=1,\ldots,\card{S_2}$. Thus, $\mbth$ is feasible for $\MP^s$ and $\tau^s\leq f(\mbth) \leq \tau^s_r(\mb{y}^1,\by^2) + \eps$. Combining with results from Proposition \ref{prop:relaxation}, we have
$$
\tau^s_r(\mb{y}^1,\by^2)\leq \tau^s\leq \tau^s_r(\mb{y}^1,\by^2) + \eps .
$$
As $\eps>0$ was arbitrarily chosen, we obtain  $\tau^s_r(\mb{y}^1,\by^2)=\tau^s$.
\end{pf}
\begin{remark}
In the case when $\card{S_i} = \binom{n+r}{r}$, the assumption that $S_i$ is $r$-generic (in the sense of Definition~\ref{def:rgeneric}) holds. Indeed, the points of $S_i$ in general position impose $\binom{n+r}{r}$ independent linear conditions, which is the maximal number of coefficients of a polynomial of degree $r$.
\end{remark}
If we select $r_0 = \min\left\{r\in \N: \binom{n+r}{r}\geq \max\{\card {S_1},\card{S_2}\}\right\}$, the condition $\card{S_i} = \binom{n+r_0}{r_0}$, $i=1,2$, does not hold in general, thus Proposition~\ref{prop:exactconvergence} cannot be directly applied.
However, we can apply the following perturbation algorithm to the initial datasets $S_1$ and $S_2$ to ensure that the rank condition generically holds:

{\bf Perturbation Algorithm}
\begin{itemize}
\item[{\bf 1.}] For $i=1,2$, replicate $(r_0 - \card{S_i})$ times an arbitrary point of  $S_i$  to obtain a new dataset $S_i'$ with $|S_i'| = \binom{n+r}{r}$.
\item[{\bf 2.}] Choose an arbitrary small $\epsilon > 0$, fixed. For $i=1,2$ and  each $x \in S_i'$, generate a random unit vector $\tilde{u}\in \mathbb{S}^{n-1}$ from the 
rotation-invariant probability distribution on $\mathbb{S}^{n-1}$ and replace $\mb{x}$ with
$\mb{x}+ \epsilon \tilde{u}\in B(\mb{x},\epsilon)$ (where $B(\mb{x};\epsilon)$ is the ball 
centered at $\mb{x}$ and with radius $\epsilon$).
The perturbed dataset $\tilde{S}_i$ is the set of all randomly generated vectors $\tilde{\mb{x}}$.
\item[{\bf 3.}] Output $\tilde{S}_1$ and  $\tilde{S}_2$.
\end{itemize}
After applying this algorithm, the rank condition generically holds and one can apply Proposition~\ref{prop:exactconvergence} to the perturbed datasets $\tilde{S}_1$ and  $\tilde{S}_2$. 

\if{
For $i=1,2$, let us denote by 
$\tilde{\mu}^i$ a probability measure generated on $\tilde{S}_i$ after perturbation of the probability measure $\mu^i$. We also note $\tilde{\mb{y}}^i$ the moment sequence of $\tilde{\mu}^i$, for $i=1,2$. In a similar way, we define the perturbed problems $\tilde{\MP}^s$ and $\tilde{\MP}_r^s(\tilde{\mb{y}}^1,\tilde{\mb{y}}^2)$ with respective solutions $\tilde{\tau}^s$, $\tilde{\tau}^s_r(\tilde{\mb{y}}^1,\tilde{\mb{y}}^2)$.



Although the result of Proposition~\ref{prop:exactconvergence} is interesting, it is not very
useful for practical algorithms.  Problem $\tilde{\MP}_{r_0}^s(\tilde{\mb{y}}^1,\tilde{\mb{y}}^2)$ has
only two LMI constraints but its matrix size is at least $\max\{\card {S_1},\card{S_2}\}$,
which could be very large. It means that Problem $\tilde{\MP}_r^s(\tilde{\mb{y}}^1,\tilde{\mb{y}}^2)$
is still computationally difficult to solve,
 when $\card {S_1}$ or $\card {S_2}$ is large. In the next section, we will use results from this proposition to show that after generic dataset perturbations, we can find appropriate probability measures $\tilde{\mu}^i$ such that $\tilde{\tau}^s_r(\tilde{\mb{y}}^1,\tilde{\mb{y}}^2)=\tilde{\tau}^s$ for $r$ as small as $d$, which is the degree of the polynomial $\theta$.
}\fi 
 
Although the result of Proposition~\ref{prop:exactconvergence} is interesting, it is not very
useful for practical algorithms.  Problem ${\MP}_{r_0}^s({\mb{y}}^1,{\mb{y}}^2)$ has
only two LMI constraints but its matrix size is at least $\max\{\card {S_1},\card{S_2}\}$,
which could be very large. It means that Problem ${\MP}_r^s({\mb{y}}^1,{\mb{y}}^2)$
is still computationally difficult to solve,
 when $\card {S_1}$ or $\card {S_2}$ is large. However in the next section, we use 
  Proposition \ref{prop:exactconvergence} and  show that 
 we can generically find appropriate probability measures ${\mu}^i$ such that ${\tau}^s_r({\mb{y}}^1,{\mb{y}}^2)={\tau}^s$ for $r$ as small as $d$, which is the degree of the polynomial $\theta$. 
 
\subsection{Optimal Measure}
\label{ssec:optimalmeasure}
The probability measure $\mu^i$ is defined in \refs{eq:mudef} as $\mu^i =
 \sum_{\mb{x}\in S_i}\mu^i_{\mb{x}}\delta_{\mb{x}}$ with $\sum_{\mb{x}\in
 S_i}\mu^i_{\mb{x}}=1$ and $\mu^i_{\mb{x}}\geq 0$ for all $\mb{x}\in S_i$.
 Let $\mbmu^i = \left(\mu^i_{\mb{x}_1},\ldots,\mu^i_{\mb{x}_{\card S_i}}\right)$, we have, $\mbmu^i \in M_{\card {S_i}}$, where $M_{\card {S_i}} = \left\{\mb{x}\in \R^{\card {S_i}}_+:\sum_{i=1}^{\card {S_i}}x_i=1\right\}$.\footnote{We use bold notation $\mbmu^i$ for the weight vector of the measure $\mu^i$, in adequation with the bold notation $\mbth$ for the coefficient vector of the polynomial $\theta$.} Each probability measure $\mu^i$ can then be represented equivalently by a vector $\mbmu^i \in M_{\card {S_i}}$. Thus, the optimal value of  Problem ${\MP}^s_r(\mb{y}^1,\by^2)$ can also be expressed as $\tau^s_r(\mbmu^1,\mbmu^2)\equiv \tau^s_r(\mb{y}^1,\by^2)$.

Clearly, we can form infinitely many moment relaxations from different probability measures generated on $S_i$ as above. But the question is then: {\em which pair of probability measures yields the best relaxation?} With $r$ fixed, consider the following problem:
\begin{equation}
\label{eq:bestmeasure}
{\MP}^s_r \quad \left[
\rho^s_r = \sup_{\mbmu^i \in M_{\card {S_i}}}\tau^s_r(\mbmu^1,\mbmu^2)
\right] \,,
\end{equation}  
where $\rho^s_r$ is the optimal value of ${\MP}^s_r$.
We then immediately have the following result
\begin{proposition}
\label{lma:bestmeasure}
Let $\MP^s$ and ${\MP}^s_r$ be as in \refs{eq:sprob} and \refs{eq:bestmeasure} respectively. Then
$
\rho^s_r \leq \tau^s.
$
\end{proposition}
\begin{pf}
From Proposition \ref{prop:relaxation}, we have $\tau^s_r(\mbmu^1,\mbmu^2)\equiv \tau^s_r(\mb{y}^1,\by^2)\leq \tau^s$ for all $\mbmu^i\in M_{\card {S_i}}$, $i=1,2$. 
Therefore $\sup_{\mbmu^i\in M_{\card {S_i}}}\tau^s_r(\mbmu^1,\mbmu^2)\leq \tau^s$.

\end{pf}

We are interested in finding the minimum value of the moment order $r$ that turns  the above inequality into an equality. We observe that 
in view of the convexity of $f$ and $\Theta$, the optimal solution of $\MP^s$ depends only on some (possibly small) subsets of $S_1 \subset \R^n$ and $S_2 \subset \R^n$.  Indeed, the following result was 
proved by Ben-Tal et al. \cite{Ben-Tal79}: 

\begin{theorem}{
[Ben-Tal et al. \cite[Theorem 3.1]{Ben-Tal79}]}
\label{thr:extremum}
Consider the problem
$$
{\MP} \quad \left[
\ba {rl}
\inf & f(\mbth)\\
\st & g_k(\mbth,\bx)\leq 0,\quad \bx\in S_k,\quad k=1,\ldots,m,\\
\quad & \mbth\in\Theta,
\ea
\right]
$$
and assume that
\begin{enumerate}
\item[(A1)] the set $\Theta \subset \R^t$ is convex with non empty interior,
\item[(A2)] the function $f$ is continuous and convex on $\Theta$,
\item[(A3)] the function $g_k$ is continuous in $\bx$,
\item[(A4)] for all $k \in \N$, the function $g_k$ is continuous and convex in $\mbth$ on $\Theta$ and the set $\{\mbth: g_k(\mbth,\bx)<0\}$ is open, for each $\bx\in S_k$,  
\item[(A5)] (Slater condition) the set $\{\mbth\in\R^t: g_k(\mbth,\bx)<0,\bx\in S_k,k=1,\ldots,m\}$ is nonempty.
\end{enumerate}
Let $\mbth^*$ be a feasible solution of $\MP$,  $S_k(\mbth^*)=\{\bx\in S_k: g_k(\mbth^*,\bx)=0\}$, and $K^*=\{k:S_k(\mbth^*)\neq\emptyset\}$. Then $\mbth^*$ is an optimal solution of $\MP$ if and only if there is a set $S^*\subset \cup_{k\in K^*} S_k(\mbth^*)$ with at most $t$ elements such that $\mbth^*$ is the optimal solution of the problem:
$$
{\MP}^* \quad \left[
\ba {rl}
\inf & f(\mbth)\\
\st & g_k(\mbth,\bx)\leq 0,\quad \bx\in S^* \cap S_k  ,\quad k\in K^*,\\
\quad & \mbth\in\Theta.
\ea
\right]
$$
\end{theorem}
Using Theorem \ref{thr:extremum} with $g_1 (\mbth,\bx) := -\mbth(\bx) $ and $g_2 (\mbth,\bx) := \mbth(\bx) $, we see that our initial problem of separating the two datasets $S_1$ and $S_2$ boils down to separating the two datasets $S^* \cap S_1$ and $S^* \cap S_2$, of smaller size, bounded by $\binom{n+d}{n}$. Our aim is then to apply Proposition~\ref{prop:exactconvergence} in order to solve exactly this equivalent problem. To do so, we need the two initial datasets $S_1$ and $S_2$ to fulfill genericity conditions, which are slightly stronger than the one stated in Definition~\ref{def:rgeneric}.
\if{
\begin{definition}[$d$-genericity]
\label{def:dgeneric}
We say that $S$ is said ``$d$-generic'' when each subset of $S$ of size less than $\binom{n+d}{n}$ does not belong to the level set of any polynomial of degree at most $d$.
\end{definition}
}\fi
\begin{theorem}
\label{thr:bestorder}
Let $\MP^s$, ${\MP}^s_r$ be defined as in \refs{eq:sprob}, \refs{eq:bestmeasure} respectively, whose variables are the coefficients of a degree $d$ polynomial $\theta \in \Theta$. 
Assume that $\Theta$ is convex, $f$ is convex on $\Theta$, Slater condition is satisfied and that each subset of $S_i$ of size less than $\binom{n+d}{n}$ for $i=1,2$ is $d$-generic (in the sense of Definition~\ref{def:rgeneric}). If ${\MP}^s$ is solvable, then 
the following generically holds for all $r\geq d$:
$$
{\rho}^s_r = {\tau}^s .
$$
\end{theorem}
\begin{pf}
Let us consider the separation problem ${\MP}^s$.  
%
In our context $g_k(\theta,\mb{x})=\sum_{\alpha}\theta_\alpha \mb{x}^\alpha$ and therefore $g_k$ satisfies (A3) and (A4). 
As ${\MP}^s$ is solvable with optimal solution $\mbth^*$, we can apply the results of Theorem~\ref{thr:extremum} with $t = \binom{n+d}{n}$. Thus, there exists a set $S^*$ such that $\mbth^*$ is the optimal solution of the reduced problem ${\MP}^{*}$ associated to ${\MP}^s$:
$$
{\MP}^* \quad \left[
 \ba {rl}
\inf & f(\mbth)\\
\st & \theta(\bx)\geq 0,\quad \bx\in S^* \cap S_1 ,\\
\quad & \theta(\bx) \leq 0,\quad \bx\in S^* \cap S_2 ,\\
\quad & \mbth\in\Theta.
\ea
\right]
$$
In general, the cardinal $c_i$ of $S^* \cap S_i$ will be strictly less than $t = \binom{n+d}{n}$ and we cannot apply Proposition~\ref{prop:exactconvergence} to the reduced problem ${\MP}^{*}$. However, we can modify ${\MP}^{*}$ by considering two sets $S_i^*$ of cardinal $t$, obtained after adding $t - c_i$ points of $S_i$ to $S^* \cap S_i$, for $i=1,2$. 
Since $S^* \cap S_i \subseteq S_i^*$ for $i=1,2$, if $\theta$ (resp.~$-\theta$) is nonnegative on $S_1^*$ (resp.~$S_2^*$) then it is also nonnegative on $S^* \cap S_1$ (resp.~$S^* \cap S_2$). In other words, it is sufficient to separate $S_1^*$ and $S_2^*$ with the level set of $\theta$ since the same level set separates $S^* \cap S_1$ and $S^* \cap S_2$ (and thus $S_1$ et $S_2$).

This leads to the following problem:

$$
{\MP}^s(S_1^*,S_2^*) \quad \left[
 \ba {lrl}
{\tau}^s(S_1^*,S_2^*) = & \inf & f(\mbth)\\
& \st & \theta(\bx)\geq 0,\quad \bx\in S_1^* ,\\
& \quad & \theta(\bx) \leq 0,\quad \bx\in S_2^* ,\\
& \quad & \mbth\in\Theta.
\ea
\right]
$$
where ${\tau}^s(S_1^*,S_2^*)$ denotes the optimal value of ${\MP}^s(S_1^*,S_2^*)$. 
In other words, the problem of separating $S_1$ and $S_2$ is equivalent to the problem of separating two datasets $S_1^*$ and $S_2^*$ of smaller size  (but of course, $S_1^*$ and $S_2^*$ are not known in advance).

Note that ${\tau}^s(S_1^*,S_2^*)$ is equal to the optimal value of ${\MP}^{*}$, which is also equal to the optimal value of ${\MP}^s$. Both optimal values are reached at the same $\mbth^*$ and $\tau^s = \tau^s(S_1^*,S_2^*)$.

Then we choose probability measures ${\mu}^i_0$ (with the moment vector ${\mb{y}}^i_0$) supported exactly on the whole set $S_i^*$, that is, for all $\mb{x} \in S_i$, $\mu^i_{\mb{x}}> 0$ if and only if $\mb{x}\in S_i^*$, for $i=1,2$. 
Clearly, $\mbmu^i_0\in M_{\card {{S}_i}}$, for $i=1,2$, thus ${\tau}^s_d({\mbmu}^1_0,{\mbmu}^2_0)\leq \rho^s_d$.

The probability measure ${\mu}^i_0$ is supported on the whole set $S_i^*$, for $i=1,2$, thus $\MP^s_d({\mb{y}}^1_0,{\mb{y}}^2_0)$ is also a moment relaxation  of Problem ${\MP}^s(S_1^*,S_2^*)$.
Thanks to the $d$-genericity assumptions on all subsets of $S_1$ and $S_2$, the two datasets $S_1^*$ and $S_2^*$ do not belong to the level set of any polynomial of degree at most $d$.
Therefore, we can apply the result from Proposition \ref{prop:exactconvergence}, yielding ${\tau}^s_d({\mb{y}}^1_0,{\mb{y}}^2_0) \equiv \tau^s_d({\mbmu}^1_0,{\mbmu}^2_0) = {\tau}^s(S_1^*,S_2^*)$ since the points of $S_i^*$ are in generic position and $|S_i^*| = t = \binom{n+d}{n}$, for $i=1,2$.

From Proposition \ref{lma:bestmeasure}, we have $ {\rho}^s_d \leq {\rho}^s_r\leq {\tau}^s$ for all $r \geq d$.
 From these inequalities and equalities, we have: ${\tau}^s={\tau}^s(S_1^*,S_2^*)=
{\tau}^s_d ({\mbmu}^1_0,{\mbmu}^2_0)\leq {\rho}^s_d \leq {\rho}^s_r \leq {\tau}^s$, for all $r\geq d$. 

\if{
The polynomial $\theta(\bx)$  is linear w.r.t.~the coefficient vector $\mbth$. If ${\MP}^s$ is solvable, then one can perturbate the datasets $S_1$ and $S_2$ using the perturbation algorithm so that:
\begin{itemize}
\item Problem $\tilde{\MP}^s$ is also solvable,
\item the result from Theorem~\ref{thr:extremum} can be applied with the set $\tilde{S}^*$ as for  the reduced problem $\tilde{\MP}^{*}$ associated to $\tilde{\MP}^s$,
\item one has $\card{\tilde{S}^*\cap \tilde{S}_i} = t = \binom{n+d}{d}$, for $i=1,2$.
\end{itemize}
}\fi
\end{pf}

Note that after running the perturbation algorithm from~\S~\ref{ssec:exactconvergence}, the $d$-genericity assumption of Theorem~\ref{thr:bestorder} is fulfilled for all subsets of $S_1$ and $S_2$.
The result is that, with $r$ as small as $d$, some moment relaxation ${\MP}^s_r({\mb{y}}^1,{\mb{y}}^2)$ is equivalent to Problem ${\MP}^s$, given that the appropriate probability measures $\mu^1$ and $\mu^2$ are used. These appropriate measures are uniformly supported on two datasets $S_1^*$ and $S_2^*$ of smaller size. If we would know these smaller datasets, we could easily separate $S_1$ and $S_2$ by solving the equivalent reduced problem.
Notice that for instance with $d=4$ and $n=3$ ($2$ clouds of 3D-points), one is left with 
$35\times 35$ psd matrices only (with $2D$-points then the size of each matrix drops to $15\times 15$).

The goal of the next section is to propose a practical iterative algorithm to compute separating candidate polynomials.
%

\section{Practical Algorithm}
\label{sec:algorithms}
\subsection{Algorithm}
\label{ssec:algorithm}
The key question is how to select the optimal probability measures for Problem ${\MP}^s_r$. The proof of Theorem \ref{thr:bestorder} suggests that in order to  find the optimal probability measures, we need to find  the set of points $S_i^*$, for $i=1,2$, that defines the optimal solution of Problem ${\MP}^s$. 
We propose a practical iterative algorithm to select the optimal probability measures.
At step $k$, one provides subsets $S_1^k \subseteq S_1$ and $S_2^k \subseteq S_2$ of size potentially larger than $\binom{n+d}{n}$, and defines two probability measures $\mu^1$ and $\mu^2$, supported on $S_1^k$ and $S_2^k$ respectively. Then one solves Problem ${\MP}^s_r({\mb{y}}^1,{\mb{y}}^2)$ to obtain a candidate polynomial  $\theta_k$ to separate $S_1$ and $S_2$.
By Theorem~\ref{thr:bestorder}, theoretically it would be enough 
to consider subsets $S_1^k$ and $S_2^k$ of size exactly equal to $\binom{n+d}{n} = |S_1^*| = |S_2^*|$. But for practical efficiency, we may and will tolerate subsets $S_1^k$ and $S_2^k$ of size potentially larger than $\binom{n+d}{n}$.
If the algorithm terminates at iteration $k = K$, then $S_1^K \supseteq S_1^*$ and $S_2^K \supseteq S_2^*$, for $S_1^*$ and $S_2^*$ as in Theorem~\ref{thr:bestorder}, and the resulting $\theta_K$ separates $S_1$ and $S_2$.

In each iteration, we solve Problem ${\MP}^s_r({\mb{y}}^1,{\mb{y}}^2)$ with different ${\mu}^1$ and ${\mu}^2$ until we (possibly) find the optimal probability measures.

\noindent
The main algorithm is described as follows:

{\bf Main Algorithm}
\begin{itemize}
\item[{\bf 1.}] Initialization: 
set $k\leftarrow 0$, ${S}^k \leftarrow {S_1} \cup {S_2}$, $r\leftarrow d$.
\item[{\bf 2.}] Create ${\mu}^i_k$ uniformly over ${S}^k\cap {S}_i$. Solve ${\MP}^s_r({\mb{y}}^1,{\mb{y}}^2)$. Obtain optimal solution $\theta_k$.
\item[{\bf 3.}] Form the set of outside points $O^k:=\left\{\mb{x}\in {S}_1: \theta_k(\mb{x})< 0\right\}\cup \left\{\mb{x}\in {S}_2: \theta_k(\mb{x})> 0\right\}$. If $O^k=\emptyset$, STOP. Return $\theta = \theta_k$. 
\item[{\bf 4.}] Update: $k\leftarrow k+1$, ${S}^k\leftarrow \left\{\mb{x}\in {S}_1: \theta_k(\mb{x})\leq 0\right\}\cup \left\{\mb{x}\in {S}_2: \theta_k(\mb{x})\geq 0\right\}$. Go to step \textbf{2}.
\end{itemize}
The update rule for supporting sets is based on the fact that the set of points outside the current optimal set, obtained from the moment relaxation, is likely to contain points that define the optimal separation (or covering) set. This is also the reason why ${S}^0$ is selected as ${S}_1\cup {S}_2$, which helps to separate critical and non-critical points right after the first iteration. After a supporting set is created, all points in the set are to be equally considered; therefore, uniform probability measures are used to form the moment relaxation in each iteration. 
\begin{proposition}
\label{prop:algo}
Let us assume that the main algorithm terminates. Then, we obtain an optimal solution of Problem $\MP^s$.
\end{proposition}
\begin{pf}
Suppose that the algorithm terminates at iteration $K$.
Then, the set of outside
points, $O^K= \left\{\mb{x}\in {S}_1: \theta_K(\mb{x})< 0\right\} \cup \left\{\mb{x}\in {S}_2: \theta_K(\mb{x})> 0\right\}$, is empty. Let us define ${S}^K_i := {S}^K \cap {S}_i$, for $i=1,2$. 
Since $\theta_K$ is an optimal solution of Problem ${\MP}^s_r({\mb{y}}^1,{\by}^2)$ with the uniform distribution ${\mu}^i_K$ over ${S}^K_i \subset {S}_i$, for $i=1,2$, then $\theta_K$ is a feasible solution of Problem ${\MP}^s$.  Hence, we  have ${\tau}^s \leq  {\tau}^s_r({\mbmu}_K^1,{\mbmu}_K^2)$, the optimal value of Problem ${\MP}^s_r({\mb{y}}^1,{\by}^2)$. 
On the other hand, Problem ${\MP}^s_r({\mb{y}}^1,{\by}^2)$ is a relaxation of Problem ${\MP}^s({S}^K_1,{S}^K_2)$, the separation problem 
constructed over ${S}_i^K$ instead of ${S}_i$, $i=1,2$. 
Therefore, $\tau^s_r({\mb{y}}_K^1,{\mb{y}}_K^2)\leq {\tau}^s({S}^K_1,{S}^K_2)$, as
${\tau}^s({S}^K_1,{S}^K_2)$ is the optimal value of Problem ${\MP}^s({S}^K_1,{S}^K_2)$. 
In addition, ${S}^K_i \subset {S}_i$, for $i=1,2$, 
thus ${\tau}^s({S}^K_1,{S}^K_2)\leq {\tau}^s$. 
Combining these inequalities, we obtain:
\[ 
{\tau}^s \leq  {\tau}^s_r({\mbmu}_K^1,{\mbmu}_K^2) \leq   {\tau}^s({S}^K_1,{S}^K_2)\leq {\tau}^s \,,
\]
which proves that the
optimal solution  of Problem ${\MP}^s_r({\mb{y}}^1,{\by}^2)$ is also optimal for Problem ${\MP}^s$.
\end{pf}

Despite the fact that our algorithm is a heuristic  and has no guarantee to terminate, the results established in the following sections show that this algorithm often terminates in practice after a few iterations.

In the sequel, we consider the minimum-volume covering ellipsoids problem (Section~\ref{ssec:mvce}) and the separation problem (Section~\ref{ssec:spve}) via ellipsoids. Some computational results with this iterative algorithm for these two problems are reported.

\subsection{Minimum-Volume Covering Ellipsoid Problem}
\label{ssec:mvce}
\subsubsection{Problem Formulations}
\label{sssec:mvceform}
The minimum-volume covering ellipsoid problem involves only one dataset. Let $S \subset \R^n$ be a finite set of points, $S =\left\{\mb{x}_1,\ldots,\mb{x}_t\right\}$, where $t= \card S$. We assume that the affine hull of $\mb{x}_1,\ldots,\mb{x}_t$ spans $\R^n$, which will guarantee any ellipsoids that cover all the points in $S$ have positive volume.

The ellipsoid $\Om \in \R^n$ to be determined can be written
$$
\Om := \left\{\mb{x}\in \R^n : (\bx-\mb{d})'\mb{Q}(\bx-\mb{d}) \leq 1 \right\},
$$
where $\mb{Q}\succ 0$, $\mb{Q}=\mb{Q}'$. The volume of $\Om$ is proportional to $\det\mb{Q}^{-1/2}$; therefore, the minimum-volume covering ellipsoid problem can be formulated as a maximum determinant problem (see Vandenberghe et al.~\cite{Vandenberghe98} and recent survey by Todd~\cite{Todd16} for more details) as follows:
\begin{equation}
\label{eq:omvce}
\MP \bydef \left[
\ba{rl}
\inf_{\mb{Q},\mb{d}} & \det \mb{Q}^{-1/2}\\
\st & (\bx-\mb{d})'\mb{Q}(\bx-\mb{d}) \leq 1, \quad \forall \mb{x}\in S,\\
\quad & \mb{Q} = \mb{Q}' \succ \mb{0}.
\ea
\right]
\end{equation}
Let $\mb{A}=\mb{Q}^{1/2}$ and $\mb{a}=\mb{Q}^{1/2}\mb{d}$. Then $\MP$ is equivalent to a convex optimization problem with $\card S + 1$ LMI constraints in the unknown variables $\mb{A}$ and $\mb{a}$. Indeed, each constraint $(\bx-\mb{d})'\mb{Q}(\bx-\mb{d}) \leq 1$ for $\mb{x}\in S$ can be rewritten as follows:
$$
\left[
\ba{cc}
\mb{I} & \mb{A}\mb{x}-\mb{a}\\
(\mb{A}\mb{x}-\mb{a})' & 1
\ea
\right]\succeq \mb{0}, \quad \mb{x}\in S,
$$
where $\mb{I}$ is the $n\times n$ identity matrix.

Instead of using $\mb{A}$ and $\mb{a}$, let consider $\mb{b} = 2\mb{Q}\mb{d}$ and $\mb{c} = 1 - \frac{1}{4}\mb{b}'\mb{Q}^{-1}\mb{b}$, then the ellipsoid $\Om$ can be written:
$$
\Om := \left\{\mb{x}\in \R^n : -\mb{x}'\mb{Q}\mb{x}+\mb{b}'\mb{x} + c \geq 0 \right\}.
$$
The minimum-volume covering ellipsoid problem is formulated as follows:
\begin{equation}
\label{eq:lmvce}
\ba{rl}
\inf_{\mb{Q},\mb{b},c} & \det \mb{Q}^{-1/2}\\
\st & -\mb{x}'\mb{Q}\mb{x}+\mb{b}'\mb{x} + c \geq 0, \quad \forall \mb{x}\in S,\\
\quad & \mb{c} = 1 - \frac{1}{4}\mb{b}'\mb{Q}^{-1}\mb{b},\\
\quad & \mb{Q} = \mb{Q}' \succ \mb{0}.
\ea
\end{equation}
Consider the relaxation 
\begin{equation}
\label{eq:rlmvce}
\ba{rl}
\inf_{\mb{Q},\mb{b},c} & \det \mb{Q}^{-1/2}\\
\st & -\mb{x}'\mb{Q}\mb{x}+\mb{b}'\mb{x} + c \geq 0, \quad \forall \mb{x}\in S,\\
\quad & \mb{c} \leq 1 - \frac{1}{4}\mb{b}'\mb{Q}^{-1}\mb{b},\\
\quad & \mb{Q} = \mb{Q}' \succ \mb{0}.
\ea
\end{equation}
\begin{lemma}
\label{lem:vol}
Any optimal solution $(Q^*,b^*,c^*)$ of Problem (\ref{eq:rlmvce}) is an optimal solution of Problem (\ref{eq:lmvce}). 
\end{lemma}
\begin{pf}
Since Problem (\ref{eq:rlmvce}) is an relaxation of Problem (\ref{eq:lmvce}), we just need to prove that $(Q^*,b^*,c^*)$ is a feasible solution of Problem (\ref{eq:lmvce}). 

Suppose there exists an optimal solution $(Q,b,c)$ of Problem (\ref{eq:rlmvce}) that satisfies the inequality $\gamma = c + \frac{1}{4}\mb{b}'\mb{Q}^{-1}\mb{b}<1$. We have:
$$
-\mb{x}'\mb{Q}\mb{x}+\mb{b}'\mb{x} + c \geq 0 \Leftrightarrow c + \frac{1}{4}\mb{b}'\mb{Q}^{-1}\mb{b} \geq (\mb{Q}^{1/2}\bx - \frac{1}{2}\mb{Q}^{-1/2}\mb{b})'(\mb{Q}^{1/2}\bx - \frac{1}{2}\mb{Q}^{-1/2}\mb{b})\geq 0.
$$
If we assume that $\card{S}>1$ then we have $\gamma = c + \frac{1}{4}\mb{b}'\mb{Q}^{-1}\mb{b}> 0$. Thus $0<\gamma<1$.

Let us consider the solution $(\tilde{\mb{Q}},\tilde{\mb{b}},\tilde{c})$ that satisfies $\mb{Q}=\gamma\tilde{\mb{Q}}$, $\mb{b}=\gamma\tilde{\mb{b}}$, and $\tilde{c}=1-\frac{1}{4}\tilde{\mb{b}}'\tilde{\mb{Q}}^{-1}\tilde{\mb{b}}$, we have:
$$
(\mb{Q}^{1/2}\bx - \frac{1}{2}\mb{Q}^{-1/2}\mb{b})'(\mb{Q}^{1/2}\bx - \frac{1}{2}\mb{Q}^{-1/2}\mb{b})= \gamma (\tilde{\mb{Q}}^{1/2}\bx - \frac{1}{2}\tilde{\mb{Q}}^{-1/2}\tilde{\mb{b}})'(\tilde{\mb{Q}}^{1/2}\bx - \frac{1}{2}\tilde{\mb{Q}}^{-1/2}\tilde{\mb{b}})
$$ 
Thus
$$
-\mb{x}'\mb{Q}\mb{x}+\mb{b}'\mb{x} + c \geq 0 \Leftrightarrow (\tilde{\mb{Q}}^{1/2}\bx - \frac{1}{2}\tilde{\mb{Q}}^{-1/2}\tilde{\mb{b}})'(\tilde{\mb{Q}}^{1/2}\bx - \frac{1}{2}\tilde{\mb{Q}}^{-1/2}\tilde{\mb{b}}) \leq 1,
$$
or we have $-\mb{x}'\tilde{\mb{Q}}\mb{x}+\tilde{\mb{b}}'\mb{x} + \tilde{c} \geq 0$ for all $\bx\in S$. Therefore, the solution $(\tilde{\mb{Q}},\tilde{\mb{b}},\tilde{c})$ is a feasible for Problem (\ref{eq:rlmvce}). However, we have:
$$
\mb{Q} = \gamma\tilde{\mb{Q}} \Rightarrow \det\tilde{\mb{Q}}^{-1/2} = \gamma^{n/2}\det\mb{Q}^{-1/2} <\det\mb{Q}^{-1/2}. 
$$
This contradicts the fact that $(\mb{Q},\mb{b},c)$ is an optimal solution of Problem (\ref{eq:rlmvce}). Thus we must have $\mb{c} = 1 - \frac{1}{4}\mb{b}'\mb{Q}^{-1}\mb{b}$ or $(\mb{Q},\mb{b},c)$ is a feasible (optimal) solution of Problem (\ref{eq:lmvce}).
\end{pf}

Using Lemma \ref{lem:vol} and the following fact:
$$
\left\{
\ba{l}
\frac{1}{4}\mb{b}'\mb{Q}^{-1}\mb{b} \leq 1-c\\
\mb{Q}=\mb{Q}\succeq \mb{0}
\ea
\right. \Leftrightarrow
\left[
\ba{cc}
\mb{Q} & \frac{1}{2}\mb{b}\\
\frac{1}{2}\mb{b}' & 1-c
\ea
\right] \succeq \mb{0},
$$
we can then formulate the minimum-volume covering problem as the following maximum determinant problem with $\card{S}$ linear constraints:
\begin{equation}
\label{eq:mvce}
\MP \bydef \left[
\ba{rl}
\inf_{\mb{Q},\mb{b},c} & \log \det \mb{Q}^{-1}\\
\st & -\mb{x}'\mb{Q}\mb{x}+\mb{b}'\mb{x} + c \geq 0, \quad \forall \mb{x}\in S,\\
\quad & \left[
\ba{cc}
\mb{Q} & \frac{1}{2}\mb{b}\\
\frac{1}{2}\mb{b}' & 1-c
\ea
\right] \succeq \mb{0}.
\ea
\right]
\end{equation}
Clearly, with this formulation, the minimum-volume covering ellipsoid problem is one of the covering problems $\MP^c$ as shown in \refs{eq:cprob}.

\subsubsection{Computational Results}
\label{sssec:mvceresults}
We have  implemented the algorithm presented in Section \ref{ssec:algorithm} with $r=2$ for the minimum-volume covering ellipsoid problem (we just need to set $S_2=\emptyset$). Datasets are generated using several independent normal distributions to represent data from one or more clusters. The data are then affinely transformed so that the geometric mean is the origin and all data points are in the unit ball. This affine transformation is done to make sure that data samples have the same magnitude. Computation is done in Matlab 8.5.0.197613 (R15a, SP3) with general-purpose YALMIP 3 interface~\cite{YALMIP} and SDPT3 3.4.0 solver~\cite{sdpt3} on on an Intel Core i7-5600U CPU ($2.60\, $GHz). Clearly, this algorithm can be implemented with SeDuMi~\cite{sedumi} or \emph{maxdet} solver~\cite{Vandenberghe98} in particular for this determinant maximization problem.
We have also implemented a variant of the minimum-volume covering ellipsoid to obtain level sets of quartic polynomials. This variant is obtained by replacing $\mb{x} = (x_1, \dots, x_n)$ with the vector of degree-two monomials $v_1(\mb{x}):=(x_1^2, x_1 x_2, x_2^2, \dots, x_n^2)$ in Problem~\eqref{eq:mvce}. That is,  the function $f(\mbth)$ to minimize reads $\log\,{\rm det}\,\mb{Q}^{-1}$ with $\mb{Q}\succeq0$
and with:
\[x\mapsto \theta(\mb{x})\,:= -\,v_1(\mb{x})^T\mb{Q}v_1(\mb{x}) +\mb{b}^T v_1(\mb{x})+c.\]
for some matrices $\mb{Q}$, vector $\mb{b}$, and scalar $c$. We obtain very similar results after choosing $\theta(\mb{x})= 1 -w_1 (\mb{x})'\mb{Q} w_1(\mb{x})$, with $w_1(\mb{x})$ being the vector of all monomials with degree at most $2$.

\begin{figure}[!htbp]
\centering
\subfigure[$d=2$]{
\includegraphics[scale=\sizesmallfig]{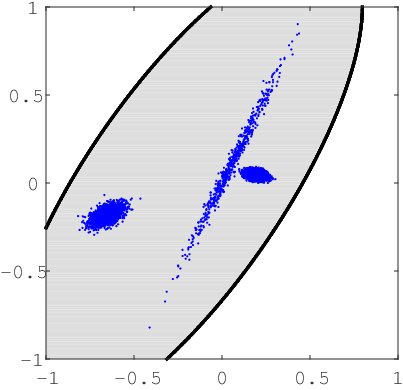}}
\subfigure[$d=4$]{
\includegraphics[scale=\sizesmallfig]{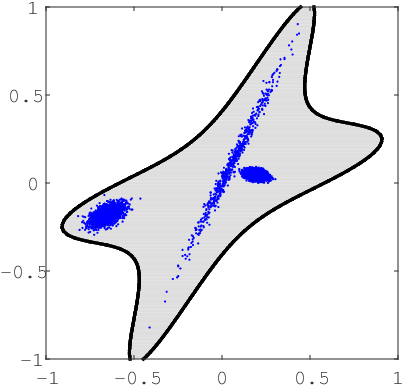}}
\caption{Minimum-volume covering ellipsoid ($d=2$) and quartic ($d=4$) for a $10000$-point dataset.}
\label{fig:mvce}
\end{figure}

The test cases show that the algorithm works well with data in two or three dimensions. Figure \ref{fig:mvce} shows the minimum-volume covering ellipsoids and quartic for a $10000$-point dataset on the plane. This figure also indicates that when the degree $d$ of the polynomial $\theta$ increases, the corresponding level set provides tighter approximation of the dataset.
When $n=3$, we have run the algorithm for datasets with up to $100,000$ points. The number of iterations we need is about $6$ and it decreases when we decrease the number of points to be covered. We also have results for datasets with $10,000$ points when $n=10$. However, the time to prepare moment matrices increases significantly in terms of the dimension. We need to prepare $O(n^2)$ square matrices of size $O(n^2)$ as data input for the relaxation if $r=2$. If the probability measure is supported on $m$ points, then the total computational time to prepare all necessary moment matrices is proportional to $O(n^7m)$. Clearly, this algorithm is more suitable for datasets in low dimensions with a large number of points. The computational time could be reduced significantly if we implement additional heuristic to find a good initial subset instead of the whole set. A problem-specific SDP code that exploits the data structure of the relaxation could be useful for datasets in higher dimensions.

\subsection{Separation Problem via Ellipsoids}
\label{ssec:spve}
\subsubsection{Problem Formulation}
\label{sssec:spveform} 
The separation problem via ellipsoids with two datasets $S_1\subset \R^n$ and $S_2\subset \R^n$ is to find an ellipsoid that contains one set, for example, $S_1$, but not the other, which is $S_2$ in this case. If we represent the ellipsoid as the set $\Om := \left\{\mb{x}\in \R^n : -\mb{x}'\mb{Q}\mb{x}+\mb{b}'\mb{x} + c \geq 0 \right\}$ with $\mb{Q}=\mb{Q}'\succ \mb{0}$, then similar to the minimum-volume covering ellipsoid problem, we can formulate the separation problem as follows:
\begin{equation}
\label{eq:spve}
\MP \bydef \left[
\ba{rl}
\inf_{\mb{Q},\mb{b},c} & \log \det \mb{Q}^{-1}\\
\st & -\mb{x}'\mb{Q}\mb{x}+\mb{b}'\mb{x} + c \geq 0, \quad \forall \mb{x}\in S_1,\\
\quad & -\mb{x}'\mb{Q}\mb{x}+\mb{b}'\mb{x} + c \leq 0, \quad \forall \mb{x}\in S_2,\\
\quad & \left[
\ba{cc}
\mb{Q} & \frac{1}{2}\mb{b}\\
\frac{1}{2}\mb{b}' & 1-c
\ea
\right] \succeq \mb{0}.
\ea
\right]
\end{equation}
    
\subsubsection{Computational Results}
\label{sssec:spveresults}
Similar to the minimum-volume ellipsoid problem, the algorithm for this separation problem can be implemented with $r=2$. With YALMIP interface and SDPT3 solver, the logdet objective function is converted to geometric mean function, which is $-(\det\mb{Q})^{1/n}$. If the problem is feasible, the optimal solution will have $\mb{Q}\succ\mb{0}$, which means that the objective value is strictly negative. This can be considered as a sufficient condition to determine that the problem is feasible. In each iteration of the algorithm, if the optimal value is zero ($\mb{Q}=\mb{0}$, $\mb{b}=\mb{0}$, and $c=0$ is a feasible solution for the subproblem solved in each iteration), then we can stop and conclude that the problem is infeasible. Existence of the critical subset that determines the problem infeasibility can be proved using the same arguments as in the proof of Theorem \ref{thr:bestorder} for the feasibility problem:
\begin{equation}
\label{eq:fspve}
\MP \bydef \left[
\ba{rl}
\inf_{\mb{Q},\mb{b},c,d} & d\\
\st & -\mb{x}'\mb{Q}\mb{x}+\mb{b}'\mb{x} + c \geq 0, \quad \forall \mb{x}\in S_1,\\
\quad & -\mb{x}'\mb{Q}\mb{x}+\mb{b}'\mb{x} + c \leq d, \quad \forall \mb{x}\in S_2,\\
\quad & \left[
\ba{cc}
\mb{Q} & \frac{1}{2}\mb{b}\\
\frac{1}{2}\mb{b}' & 1-c
\ea
\right] \succeq \mb{0}.
\ea
\right]
\end{equation}

As for the minimum-volume covering ellipsoid problem, we have implemented a variant of the minimum-volume separation ellipsoid to separate datasets by using quartic polynomials. 
In order to test the algorithm, we generate datasets $S_1$ and $S_2$ as for the minimum-volume covering ellipsoids problem. In most cases, if we run the algorithm for $S_1$ and $S_2$, we get  infeasibility results. In order to generate separable datasets, we run the minimum-volume covering ellipsoid algorithm for $S_1$ and generate the separable set $S_2'$ from $S_2$ by selecting all points that are outside the ellipsoid. We also try to include some points that are inside the ellipsoid to test the cases when $S_1$ and $S_2'$ are separable by a different ellipsoid rather than the minimum-volume ellipsoid that covers $S_1$.

\begin{figure}[!htbp]
\centering
\subfigure[Minimum-volume covering ellipsoid]{
\includegraphics[scale=\sizesmallfig]{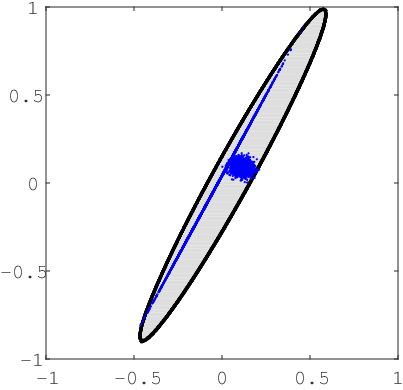}}
\subfigure[Separating ellipsoid]{
\includegraphics[scale=\sizesmallfig]{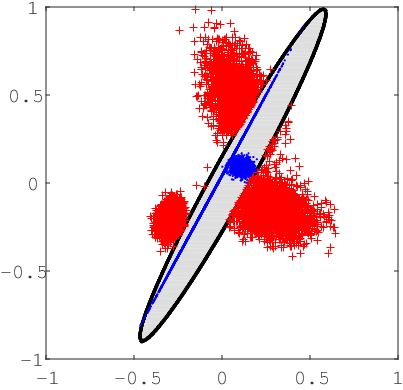}}
\caption{Separating ellipsoid is the same as the minimum-volume covering ellipsoid.}
\label{fig:mspve}
\end{figure}

\begin{figure}[!ht]
\centering
\includegraphics[scale=\sizesmallfig]{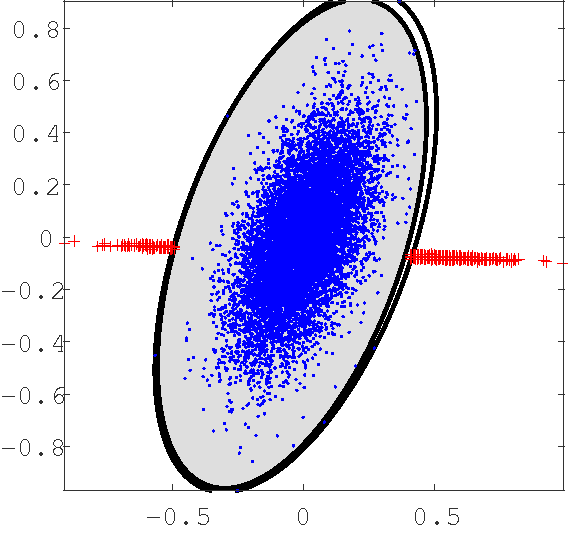}
\caption{Separating ellipsoid is different from the minimum-volume covering ellipsoid.}
\label{fig:dspve}
\end{figure}

\begin{figure}[!ht]
\centering
\subfigure[Separable datasets $S_1$ and $S_2$]{
\includegraphics[scale=\sizesmallfig]{sep0_10000_2_4.png}}\\
\subfigure[Minimum-volume ellipsoid for $S_1$]{
\includegraphics[scale=\sizesmallfig]{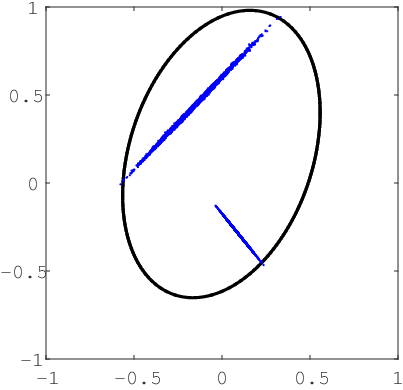}}
\subfigure[Separating quartic for $S_1$ and $S_2$]{
\includegraphics[scale=\sizesmallfig]{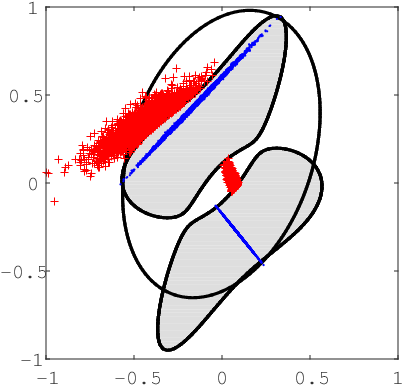}}
\caption{Separating quartic for two $10000$-point datasets.}
\label{fig:mvce4}
\end{figure}

The test cases show that the algorithm can detect problem infeasibility and in the separable case, finds an ellipsoid that separates two datasets. Figure \ref{fig:mspve} shows the separation of two datasets on the plane with $10000$ points by the minimum-volume ellipsoid, while Figure \ref{fig:dspve} represents the case when a different ellipsoid is needed to separate two particular sets. We also ran the algorithm for datasets with $n=3$ and $n=10$. Similar remarks can be made with respect to data preparation and other algorithmic issues as in Section \ref{sssec:mvceresults}. In general, the algorithm is suitable for datasets in low dimensions with a large number of points. 

Figure \ref{fig:mvce4} shows an example where there is no ellipsoid that can separate two given datasets $S_1$ and $S_2$. We indicate the border of the minimum-volume covering ellipsoid for the dataset $S_1$ as well as the level set of the separating quartic for $S_1$ and $S_2$. In such cases, one has to rely on higher degree polynomials to be able to separate the two datasets.

\renewcommand{\baselinestretch}{1.00}
\small
\bibliographystyle{plain}

\end{document}